# Convex Hull of the Quadratic Branch AC Power Flow Equations and Its Application in Radial Distribution Networks

Qifeng Li, *Member, IEEE*, and Vijay Vittal, *Fellow, IEEE*

*Abstract*—A branch flow model (BFM) is used to formulate the AC power flow in general networks. For each branch/line, the BFM contains a non-convex quadratic equality. A mathematical formulation of its convex hull is proposed, which is the tightest convex relaxation of this quadratic equation. The convex hull formulation consists of a second order cone inequality and a linear inequality within the physical bounds of power flows. The convex hull formulation is analytically proved and geometrically validated. An optimal scheduling problem of distributed energy storage (DES) in radial distribution systems with high penetration of photovoltaic resources is investigated in this paper. To capture the performance of both the battery and converter, a second-order DES model is proposed. Following the convex hull of the quadratic branch flow equation, the convex hull formulation of the nonconvex constraint in the DES model is also derived. The proposed convex hull models are used to generate a tight convex relaxation of the DES optimal scheduling (DESOS) problem. The proposed approach is tested on several radial systems. A discussion on the extension to meshed networks is provided.

## I. INTRODUCTION

THE power flow constraints are involved in many decision-making processes in power systems. For the optimization problems in distribution systems, the popular linear approximation of power flow equations, the DC power flow model [1], is no longer suitable for formulating the power flow due to the high $r/x$-ratio of the feeders. The conventional rectangular and polar AC power flow models are based on only bus variables, i.e. voltage variables. They are valid for both radial and meshed networks. However, they may introduce complexity in computation and unnecessarily high non-convexity. Some AC power flow models for radial networks have been proposed to obtain higher computational efficiency and convergence [2] and [3]. They are based on not only the bus variables but also the branch variables and, consequently, called branch flow models (BFMs) [4]. With additional angle constraints, these model can also be applied to meshed networks [4] and [5].

Non-convex constraints exist in all AC power flow models. Generally, BFMs are preferred in radial networks since the angle constraints can be neglected [6] and [7].

The convex relaxations of various AC power flow models have been studied in the context of convexifying the optimal power flow (OPF) problems. The most frequently adopted approaches are the second order cone programming (SOCP) [8], semidefinite programming (SDP) [9] relaxations and their variations [10]. When a BFM is used to model the AC power flow in radial networks, the basic SOCP relaxation dominates the basic SDP relaxation [5]. Based on the BFM in radial networks, reference [11] proposes a set of non-iterative linear cuts to obtain a tight root node for solving the OPF problem in BARON [12]; valid linear cuts are generated iteratively for tightening the SOCP relaxation of OPFs in [13]; Several convex constraints which can approximate the rank-1 constraint are introduced to obtain enhanced SDP relaxations in [7]. Based on the bus injection power flow model in polar coordinates, several linear and nonlinear cuts are introduced in [14] for strengthening the SDP relaxation in the context of OPF. As an alternative, the convex hull [15] of the non-convex constraints in the *DistFlow* model [2] is explored in this paper. The concept of convex hull is attractive since it is defined as the tightest convex relaxation of a non-convex set [16].

An essential aspect of the smart grid technologies is accommodating increasing penetration of renewable generation and storage options at both the transmission and distribution level [17]. This paper also examines the application of convex relaxations to an optimal scheduling problem of distributed energy storage (DES) in distribution systems where high penetration of photovoltaic (PV) resources are integrated. The DES in this context primarily relates to battery energy storage. The linear energy storage models considered in many related references, e.g. [6], [18], and [19], are simplified formulations of battery storage systems. To more accurately capture the performance of the DES system (including the battery and converter), a second-order model is proposed based on an equivalent circuit of the DES system.

As shown in [6], [7] and [13], the additional energy storage constraints will introduce more factors that can affect the exactness of a convex relaxation, even though these constraints are linear. Moreover, the objective function of an energy storage optimization problem is not limited to minimizing fuel cost which is a typical objective function of an OPF. It is shown in

---
Q. Li was with the Department of Electrical Engineering, Arizona State University, Tempe, AZ 85287 USA. He is now with the Department of Mechanical Engineering, Massachusetts Institute of Technology, Cambridge, MA 02139 USA (e-mail: qifengli@mit.edu).

V. Vittal is with the Department of Electrical Engineering, Arizona State University, Tempe, AZ 85287 USA (e-mail: vijay.vittal@asu.edu).

[6] and [7] that the descent direction of the objective function also impacts the exactness of the convex relaxation. With the proposed second-order DES model, the situation becomes more complicated, and requires stronger convex relaxations.

Despite years of research, it is still hard to guarantee the exact global optimal solutions of all optimization problems where the AC power flow constraints are taken into account, since some infeasible points are inevitably included in a convex relaxation of the nonconvex problem. This paper aims at proposing a tighter convex relaxation of the BFM to eliminate as many infeasible points as possible from the convex set, so that the possibility of obtaining an exact global optimal solution for actual power systems can be effectively increased. Instead of making network assumptions, the analysis in this paper is based on the original exact AC power flow model. The contributions of this paper are twofold: 1) introduce a more accurate nonlinear battery model with reactive capability considered into the optimization model of distribution systems with higher penetration of PV resources; 2) explore the convex hull of the quadratic equations in the BFM to generate a tighter convex relaxation for the above optimization problem.

The rest of the paper is organized as follow: the convex hull formulation of the non-convex quadratic constraint in the *DistFlow* model is proposed in Section II. A DESOS model with a second-order DES model is considered in Section III. In the same section, a novel convex relaxation based on the convex hull formulations is proposed for the DESOS problem. In Section IV, the proposed approach is tested in a real-world feeder as well as several IEEE standard feeders.

## II. CONVEX HULL OF AC POWER FLOW IN RADIAL NETWORKS

### A. AC Power Flow in Radial Networks

A BFM which is based on the *DistFlow* model is adopted to capture the AC power flow in general networks:

$$p_i = \sum_{k \in B_i^D} P_{ik} - \sum_{j \in B_i^U} \left( P_{ji} - r_{ij} \ell_{ji} \right) \tag{1a}$$

$$q_i = \sum_{k \in B_i^D} Q_{ik} - \sum_{j \in B_i^U} \left( Q_{ji} - x_{ij} \ell_{ji} \right) \tag{1b}$$

$$v_i = v_k + 2\left( r_{ik} P_{ik} + x_{ik} Q_{ik} \right) - \left( r_{ik}^2 + x_{ik}^2 \right) \ell_{ik} \tag{1c}$$

$$v_i \ell_{ik} = P_{ik}^2 + Q_{ik}^2 \tag{1d}$$

$$\sum_{ik \in C_l} \arctan\left( \frac{x_{ik} P_{ik} - r_{ik} Q_{ik}}{v_i - r_{ik} P_{ik} - x_{ik} Q_{ik}} \right) = 0 (\bmod\ 2\pi) \tag{1e}$$

where $i \in B$, $k \in B_i^D$. The feasible set of (1) is subject to

$$0 \leq \ell_{ik} \leq \overline{\ell}_{ik} \tag{2a}$$

$$\underline{v}_i \leq v_i \leq \overline{v}_i \tag{2b}$$

$$P_{ik}^2 + Q_{ik}^2 \leq \overline{S}_{ik}^2 \tag{2c}$$

where both $\overline{\ell}_{ik}$ and $\overline{S}_{ik}$ are related to the thermal limit of a feeder. Generally, the current limit is determined based on the conductor used in a transmission line or feeder. The relation between the thermal limit and the current limit (in pu) of a line is given as $\overline{S}_{ik} = \overline{I}_{ik} V_i^{nom} = \sqrt{\overline{\ell}_{ik} v_i^{nom}}$ ($\underline{v}_i \leq v_i^{nom} \leq \overline{v}_i$). The nomenclature for the symbols used in this section is listed in Table I.

TABLE I
NOMENCLATURE IN SECTION II

| SYMBOL | DIMENSION/QUANTITY |
|---|---|
| $B$, $B_i^D$, $B_i^U$ | Bus set of the feeder, downstream and upstream bus sets of bus $i$ respectively |
| $C_l$ | Branch set of the $l$th cycle in a meshed network |
| $L$ | Branch set of the whole system |
| $p_i$, $q_i$ | Active and reactive power injections at bus $i$ respectively |
| $P_{ik}$, $Q_{ik}$, | Active and reactive power flow in branch $ik$ respectively |
| $r_{ik}$, $x_{ik}$ | Resistance and reactance of branch $ik$ respectively |
| $v_i$, $v^{nom}$ | Square of voltage magnitude at bus $i$, and its nominal value respectively |
| $\ell_{ik}$, $\theta_{ik}$ | Square of current magnitude and angle difference in branch $ik$ respectively |
| $\overline{S}_{ik}$ | Thermal limit of line $ik$ |

Constraint (1e) can be omitted when this BFM is used to formulate the AC power flow in radial networks. For more discussions on the angle constraint (1e), please refer to [4]. As a result, the only non-convex constraints are those in (1d) and the BFM is preferred for radial networks. Moreover, there is at most one upstream bus for each bus in a radial network. Hence, each set $B_i^U$ contains only one element. A distribution system is usually operated in a radial structure under normal conditions. Consequently, this paper focuses on radial networks and the convexification of constraint (1d). A discussion on the scenario in meshed networks is offered in Subsection IV-D. Note that an analogous AC power flow model was proposed in [3] and it is equivalent to the *DistFlow* model. The proposed approach can also be applied to this model.

### B. The Convex Hull of Equation (1d)

Let $\boldsymbol{x}_{ik} = [P_{ik}\ Q_{ik}\ \ell_{ik}\ v_i]^T$ ($ik \in L$) and $\Omega_0$ denote the feasible set of equation (1d) within constraints (2), then $\Omega_0 = \{\boldsymbol{x}_{ik}\ |\ $(1d) and (2) hold.$\}$. Let ($\tilde{P}_{ik,1}$, $\tilde{Q}_{ik,1}$), ($\tilde{P}_{ik,2}$, $\tilde{Q}_{ik,2}$) be two sets of chosen values which satisfy

$$\begin{cases} \tilde{P}_{ik,1}^2 + \tilde{Q}_{ik,1}^2 = \tilde{P}_{ik,2}^2 + \tilde{Q}_{ik,2}^2 = \overline{S}_{ik}^2 \\ \left( \tilde{P}_{ik,1}, \tilde{Q}_{ik,1} \right) \neq \left( \tilde{P}_{ik,2}, \tilde{Q}_{ik,2} \right) \end{cases},$$

we then have the following proposition.

**Proposition 1.** In the $\boldsymbol{x}_{ik}$-space, the four vectors, ($\tilde{P}_{ik,1}$, $\tilde{Q}_{ik,1}$, $\overline{\ell}_{ik}$, $v_i^{nom}$), ($\tilde{P}_{ik,1}$, $\tilde{Q}_{ik,1}$, $\overline{S}_{ik}^2/\overline{v}_i$, $\overline{v}_i$), ($\tilde{P}_{ik,2}$, $\tilde{Q}_{ik,2}$, $\overline{\ell}_{ik}$, $v_i^{nom}$), and ($\tilde{P}_{ik,2}$, $\tilde{Q}_{ik,2}$, $\overline{S}_{ik}^2/\overline{v}_i$, $\overline{v}_i$), defined in the first paragraph of this subsection are linearly independent.

The proof of proposition 1 consists of showing that the determinant of the 4×4 matrix formed by taking the four vectors as its columns is non-zero. The calculation of a determinant is direct but cumbersome to show. Therefore, the proof of proposition 1 is not provided here. By substituting the four vectors into the linear equation (3) respectively it suffices to show that, in the $\boldsymbol{x}_{ik}$-space, the above four points are all located on the hyperplane which is specified by (3).

$$c_{ik}^T x_{ik} - d_{ik} = 0, \tag{3}$$

where $c_{ik} = [0\ 0\ \bar{v}_i\ \bar{\ell}_{ik}]^T$, and $d_{ik} = \bar{\ell}_{ik}(\bar{v}_i + v^{nom})$. In the $x_{ik}$-space, let $\tilde{x}_{ik}$ denote a given point in the set $\Omega_1 = \{x_{ik} \mid$ (2) and (3) hold.$\}$, and let $\tilde{x}_{ik,i}$ denote the $i$th element in vector $\tilde{x}_{ik}$, then we have the following proposition.

**Proposition 2.** For any (given) point $\tilde{x}_{ik}$, there exist two sets of values ($\tilde{P}_{ik,1}$, $\tilde{Q}_{ik,1}$), ($\tilde{P}_{ik,2}$, $\tilde{Q}_{ik,2}$) which are defined in the first paragraph of this subsection satisfying

$$\tilde{x}_{ik} = \gamma_1 \begin{bmatrix} \tilde{P}_{ik,1} \\ \tilde{Q}_{ik,1} \\ \bar{\ell}_{ik} \\ v_i^{nom} \end{bmatrix} + \gamma_2 \begin{bmatrix} \tilde{P}_{ik,1} \\ \tilde{Q}_{ik,1} \\ \bar{S}_{ik}^2/\bar{v}_i \\ \bar{v}_i \end{bmatrix} + \gamma_3 \begin{bmatrix} \tilde{P}_{ik,2} \\ \tilde{Q}_{ik,2} \\ \bar{\ell}_{ik} \\ v_i^{nom} \end{bmatrix} + \gamma_4 \begin{bmatrix} \tilde{P}_{ik,2} \\ \tilde{Q}_{ik,2} \\ \bar{S}_{ik}^2/\bar{v}_i \\ \bar{v}_i \end{bmatrix},$$

where $\gamma_i \geq 0$ and $\sum_i^4 \gamma_i = 1$ ($i = 1, \ldots, 4$).

**Proof.** Proposition 1 implies that any point in the $x_{ik}$-space can be expressed as a linear combination of the four vectors given in it. The projection of set $\Omega_1$ onto the $(\ell_{ik}\ v_i)$-plane is the line segment between ($\bar{S}_{ik}^2/\bar{v}_i$, $\bar{v}_i$) and ($\bar{\ell}_{ik}$, $v^{nom}$). It suffices to show that the projection of $\tilde{x}_{ik}$ onto the $(\ell_{ik}\ v_i)$-plane is given by

$$(\tilde{x}_{ik})_{\ell v} = \gamma_{13} \begin{bmatrix} \bar{\ell}_{ik} \\ v_i^{nom} \end{bmatrix} + \gamma_{24} \begin{bmatrix} \bar{S}_{ik}^2/\bar{v}_i \\ \bar{v}_i \end{bmatrix},$$

where $\gamma_{13}, \gamma_{24} \geq 0$ and $(\gamma_{13} + \gamma_{24} = 1)$. By carefully choosing ($\tilde{P}_{ik,1}$, $\tilde{Q}_{ik,1}$) and ($\tilde{P}_{ik,2}$, $\tilde{Q}_{ik,2}$), the projection of $\tilde{x}_{ik}$ onto the $(P_{ik}\ Q_{ik})$-plane can be described in a similar way

$$(\tilde{x}_{ik})_{PQ} = \gamma_{12} \begin{bmatrix} \tilde{P}_{ik,1} \\ \tilde{Q}_{ik,1} \end{bmatrix} + \gamma_{34} \begin{bmatrix} \tilde{P}_{ik,2} \\ \tilde{Q}_{ik,2} \end{bmatrix}$$

where $\gamma_{12}, \gamma_{34} \geq 0$ and $(\gamma_{12} + \gamma_{34} = 1)$. If $(\tilde{x}_{ik})_{PQ} = (0, 0)$, the values for $\gamma_{12}$ and $\gamma_{34}$ are unique, i.e. $\gamma_{12} = \gamma_{34} = 0.5$. The situation is pictorially described in Fig. 1.

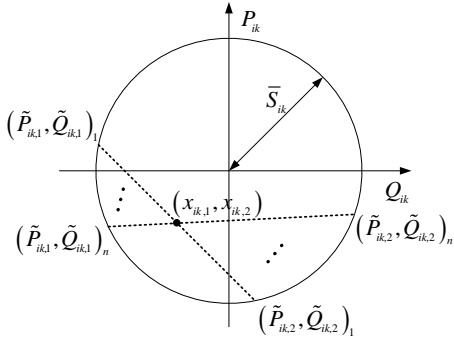

Fig. 1 The projection of of $\tilde{x}_{ik}$ onto the $(P_{ik}\ Q_{ik})$-plane.

For any point within the circle in Fig. 1, there always exist two set of values ($\tilde{P}_{ik,1}$, $\tilde{Q}_{ik,1}$) and ($\tilde{P}_{ik,2}$, $\tilde{Q}_{ik,2}$) making $\gamma_{12} = \gamma_{34} = 0.5$. Based on this finding, the proof of proposition 2 is equivalent to showing the feasibility of the following system with respect to $\gamma_i \geq 0$ ($i = 1, \ldots, 4$)

$$\begin{cases} \gamma_1 + \gamma_2 + \gamma_3 + \gamma_4 = 1 \\ \gamma_1 + \gamma_2 = \gamma_{12} = \gamma_3 + \gamma_4 = \gamma_{34} = 0.5 \\ \gamma_1 + \gamma_3 = \gamma_{13}, \gamma_2 + \gamma_4 = \gamma_{24} = 1 - \gamma_{13} \end{cases}.$$

Actually, the value of $\gamma_{13}$ is determined by the given point $\tilde{x}_{ik}$ and changes from 0 to 1. It is easy to verify that for any given $\gamma_{13}$ ($0 \leq \gamma_{13} \leq 1$), the above system is feasible. If $\gamma_{13}$ does not equal to 0 or 1, the solution is non-unique. As a result, proposition 2 is proved.

□

Note that the variables $\tilde{P}_{ik,1}$, $\tilde{Q}_{ik,1}$, $\tilde{P}_{ik,2}$, and $\tilde{Q}_{ik,2}$ mentioned in Propositions 1 and 2 are the same. They are defined right above Proposition 1. There may not exist two points that are valid for all the points in set $\Omega_1$. However, it can be guaranteed that, for any given point $\tilde{x}_{ik}$, there exist a pair of points ($\tilde{P}_{ik,1}$, $\tilde{Q}_{ik,1}$), ($\tilde{P}_{ik,2}$, $\tilde{Q}_{ik,2}$) which satisfy the conditions given in and right above proposition 1.

**Theorem.** The convex hull of $\Omega_0$ can be formulated as

$$\Omega_2 = CH(\Omega_0) = \left\{ x_{ik} \middle| \begin{array}{l} \|a^T x_{ik}\|_2 - b^T x_{ik} \leq 0 \\ c_{ik}^T x_{ik} - d_{ik} \leq 0 \\ x_{ik} \in (2) \end{array} \right\}, \tag{4}$$

where $a = [\sqrt{2}\ \sqrt{2}\ 1\ 1]^T$, $b = [0\ 0\ 1\ 1]^T$.

**Proof.** Note that, unless otherwise stated, all the discussions are within the set described by (2).

(i) $CH(\Omega_0) \subseteq \Omega_2$.
It suffices to show that

$$\{x_{ik} \mid v_i \ell_{ik} = P_{ik}^2 + Q_{ik}^2\} \subseteq \{x_{ik} \mid \|a^T x_{ik}\|_2 - b^T x_{ik} \leq 0\}.$$

Moreover,

$$\left\{ x_{ik} \middle| \begin{array}{l} v_i \ell_{ik} = P_{ik}^2 + Q_{ik}^2 \\ x_{ik} \in (2) \end{array} \right\} \subseteq \left\{ x_{ik} \middle| \begin{array}{l} v_i \ell_{ik} \leq \bar{S}_{ik}^2 \\ x_{ik} \in (2) \end{array} \right\} \subseteq \left\{ x_{ik} \middle| \begin{array}{l} c_{ik}^T x_{ik} \leq d_{ik} \\ x_{ik} \in (2) \end{array} \right\}.$$

That means $\Omega_2$ is a convex relaxation of $\Omega_0$. $CH(\Omega_0)$ is defined as the intersection of all convex relaxations of $\Omega_0$. As a result, (i) holds.

(ii) $CH(\Omega_0) \supseteq \Omega_2$.
Supposed that $a^T x_{ik} \geq \beta$ is a valid inequality for $CH(\Omega_0)$, it should be valid for all the points in $\Omega_0$. Note that "an inequality is valid for a set" means the inequality is satisfied by all the points within this set. This section tries to show that $a^T x_{ik} \geq \beta$ is also valid for all the boundaries of the convex set $\Omega_2$. Within space (2), the set $\Omega_2$ is enclosed by two sets of boundaries. They are specified by the two inequalities in (4) respectively. Within set (2), the first set of the boundary of $\Omega_2$ can be described as

$$\{x_{ik} \mid \|a^T x_{ik}\|_2 - b^T x_{ik} = 0\}, \tag{5a}$$

which is exactly $\Omega_0$. Hence, $a^T x_{ik} \geq \beta$ is valid for set (5a).
By substituting the four points into (1d), it can be observed that they all belong to $\Omega_0$. Based on proposition 2, substituting

$\tilde{x}_{ik}$ into $\alpha^T x_{ik} \geq \beta$ results in inequality (5b) which means $\alpha^T x_{ik} \geq \beta$ is valid for the second set of boundaries of $\Omega_2$, i.e. $\Omega_1$. To sum up, $\alpha^T x_{ik} \geq \beta$ is valid for all the boundary points of $\Omega_2$, which means $\alpha^T x_{ik} \geq \beta$ is also valid for the set $\Omega_2$. Hence, (ii) holds.

□

To provide an intuitive interpretation, a geometrical validation of the above theorem is given in the appendix section. Among all the existing convex relaxations, formulation (4) is the one that has been analytically proved and geometrically verified to be the convex hull of the quadratic equation (1d). A tight convex relaxation for the BFMs can be expected when (1d) is replaced with (4). As a result, the relation between the CH relaxation and the existing approaches is given in (6) where $CR(\cdot)$ denotes a general convex relaxation of a nonconvex set and $CH(\cdot) \subseteq CR(\cdot)$. Moreover, the constraints in (4) are easy to compute and do not rely on network assumptions. Therefore, we suggest the use of the CH relaxation when one needs to solve an optimization problem in radial networks accounting for the AC power flow constraints.

Generally, the extreme points of a convex hull belong to its original non-convex set. If the objective function is a convex function and monotonic over the convex hull, the optimal solution is usually located at one of the extreme points [16], implying that the optimal solution obtained by solving the convex relaxation is most likely the exact globally optimal solution of the original problem. A pictorial interpretation of the above statements is given in Fig.2. This property makes the concept of convex hull extremely attractive in the area of power system convex optimization. However, the convex hull for the feasible set of the entire branch flow model is still unknown due to the relation given in (6). Formulating the exact convex hull of the BFM given in (1a) – (1d), and (2) is the primary future work of this research.

### III. APPLICATION: OPTIMAL OPERATION OF DES

This section studies the application of the proposed approach in convexifying the DESOS problems which contain more factors that may impact the exactness of the convex re-laxation than the conventional OPF problem.

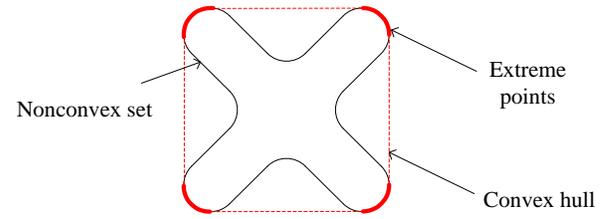

Fig. 2 A pitorial interpretation of the advatages of convex hull.

#### A. A Second Order Model of DES

In [7], a conventional storage model is used where active charge and discharge powers are defined separately and reactive power capability of the DES system is not considered. By introducing a small perturbation to the objective functions, the complementary constraints (bilinear equality constraints) on charge and discharge powers of DES units are eliminated and the resulting DES model is linear. The losses are proportional to the active power of the DES unit, which is inaccurate but acceptable.

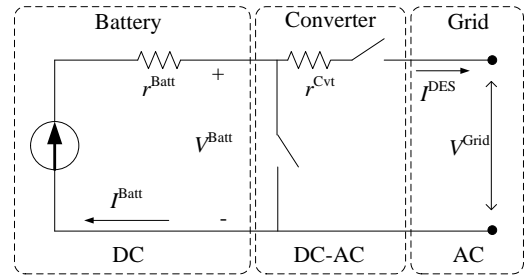

Fig. 3 The simplified equivalent circuit of a DES system.

In this paper, the reactive power capability of a DES is considered. In this case, the losses in a DES unit (including the converter and battery) are related to both active and reactive powers. To better capture the state of charge (SOC) and the loss properties of DES in this case, a second-order DES model is considered. A widely used equivalent circuit of a DES system [20] is given in Fig. 3. Based on this equivalent circuit, the

$$\alpha^T \tilde{x}_{ik} = \gamma_1 \alpha^T \begin{bmatrix} \tilde{P}_{ik,1} \\ \tilde{Q}_{ik,1} \\ \overline{\ell}_{ik} \\ v_i^{nom} \end{bmatrix} + \gamma_2 \alpha^T \begin{bmatrix} \tilde{P}_{ik,1} \\ \tilde{Q}_{ik,1} \\ \overline{S}_{ik}^2 / \overline{v}_i \\ \overline{v}_i \end{bmatrix} + \gamma_3 \alpha^T \begin{bmatrix} \tilde{P}_{ik,2} \\ \tilde{Q}_{ik,2} \\ \overline{\ell}_{ik} \\ v_i^{nom} \end{bmatrix} + \gamma_4 \alpha^T \begin{bmatrix} \tilde{P}_{ik,2} \\ \tilde{Q}_{ik,2} \\ \overline{S}_{ik}^2 / \overline{v}_i \\ \overline{v}_i \end{bmatrix} \geq \left( \sum_i^4 \gamma_i \right) \beta = \beta, \quad (5b)$$

$$\underbrace{CH \left\{ \begin{bmatrix} \text{convex constraint \#1} \\ \vdots \\ \text{convex constraint \#i} \\ \text{nonconvex constraint \#1} \\ \vdots \\ \text{nonconvex constraint \#j} \end{bmatrix} \right\}}_{\text{In future research}} \subseteq \underbrace{\left\{ \begin{matrix} \text{convex constraint \#1} \\ \vdots \\ \text{convex constraint \#1} \\ CH \begin{pmatrix} \text{nonconvex constraint \#1} \\ \vdots \\ \text{nonconvex constraint \#j} \end{pmatrix} \end{matrix} \right\}}_{\text{In this paper}} \subseteq \underbrace{\left\{ \begin{matrix} \text{convex constraint \#1} \\ \vdots \\ \text{convex constraint \#1} \\ CR \begin{pmatrix} \text{nonconvex constraint \#1} \\ \vdots \\ \text{nonconvex constraint \#j} \end{pmatrix} \end{matrix} \right\}}_{\text{In literature}} \quad (6)$$

power losses in a battery can be formulated as

$$p_i^{DESloss} = r_i^{Batt}\left(I_{i,t}^{Batt}\right)^2 + r_i^{Cvt}\left(I_{i,t}^{DES}\right)^2$$

$$= r_i^{Batt}\frac{\left(p_{i,t}^{DES}\right)^2}{v_{i,t}^{Batt}} + r_i^{Cvt}\frac{\left(p_{i,t}^{DES}\right)^2 + \left(q_{i,t}^{DES}\right)^2}{v_{i,t}}.$$

In the per unit system, $v_{i,t}^{Batt} = v_{i,t}$. Therefore, the above constraint can be rewritten as

$$p_{i,t}^{DESloss} v_{i,t} = r_i^{eq}\left(p_{i,t}^{DES}\right)^2 + r_i^{Cvt}\left(q_{i,t}^{DES}\right)^2, \quad (7)$$

where $p_{i,t}^{DES}, q_{i,t}^{DES} > 0$ means discharging, and

$$r_i^{eq} = r_i^{Batt} + r_i^{Cvt},$$

$$\left(p_{i,t}^{DES}\right)^2 + \left(q_{i,t}^{DES}\right)^2 \leq \left(\bar{S}_i^{DES}\right)^2. \quad (8a)$$

Consequently, the state of charge constraint is given as

$$\underline{E}_i^{DES} \leq E_i^{spl} - \sum_{t'=1}^{t}\left(p_{i,t'}^{DES} + p_{i,t'}^{DESloss}\right)\Delta t \leq \bar{E}_i^{DES}, \quad (8b)$$

where $i \in N_S$ and, unless otherwise stated, $t \in T$ throughout the paper. The nomenclature for the symbols used in this section is listed in Table II.

TABLE II
NOMENCLATURE IN SECTION III

| SYMBOL | DIMENSION/QUANTITY |
|---|---|
| $c_t$ | Cost of grid energy ($/MWh). In a market environment, it represents the nodal price (NP) at the distribution substation location. |
| $E_i^{spl}$ | DES energy surplus at the beginning of a day at bus $i$ |
| $\underline{E}_i^{DES}, \bar{E}_i^{DES}$ | Minimum and maximum allowed SOC of DES at bus $i$ |
| $k_i$ | Loss coefficient of the transformer at bus $i$ |
| $N_S, N_T$ | DES bus set and high voltage bus set of transformers |
| $p_{i,t}^{DES}, q_{i,t}^{DES}$ | Active and reactive outputs of DES at bus $i$ at time $t$ respectively |
| $p_{(i=1),t}^{Grid}, q_{(i=1),t}^{Grid}$ | Active and reactive grid power at hour $t$ |
| $p_{i,t}^{DESloss}$ | Power losses of DES at bus $i$ at hour $t$ |
| $p_{i,t}^L, q_{i,t}^L$ | Active and reactive load at bus $i$ at hour $t$ |
| $p_{i,t}^{PV}$ | PV generation at bus $i$ and at hour $t$ |
| $R$ | Substation MVA rating |
| $\bar{S}_i^{DES}$ | Converter MVA limit of the DES unit at bus $i$ |
| $v_i^{set}$ | Set point of the square of voltage magnitude at bus $i$ |
| $T, \Delta t$ | Operation cycle and time interval between two operations of the DES respectively |

### B. Problem Formulation

In a distribution feeder with high penetration of PV resources, the operating objectives include but are not limited to minimizing purchase cost of grid energy, network losses and voltage magnitude deviation. An operation cycle of 24 hours is considered for the DES units. The optimal charge and discharge schedule for each time interval is computed a day ahead based on 24-hour predicted load and renewable generation profiles. The topology of distribution systems is usually radial. Hence, based on the BFM given in (1) and (2), the DESOS model can be formulated as

$$\min \quad f_1 = \sum_t^T c_t p_t^{Grid} \quad (9a)$$

$$\min \quad f_2 = \sum_t^T\left(\sum_i r_{ik}\ell_{ik,t} + \sum_i^{N_T} k_i v_{i,t} + \sum_i^{N_S} p_{i,t}^{DESloss}\right) \quad (9b)$$

$$\min \quad f_3 = \sum_t^T \sum_i^B \left|v_{i,t} - v_i^{set}\right| \quad (9c)$$

s.t. (2), (7), (8) and

$$p_{i,t}^{DES} + p_{(i=1),t}^{Grid} + p_{i,t}^{PV} - p_{i,t}^L - k_{(i\in N_T)}v_{(i\in N_T),t}$$
$$= \sum_{k\in B_i^D} P_{ik,t} - \left(P_{ji,t} - r_{ij}\ell_{ji,t}\right) \quad (10a)$$

$$q_{i,t}^{DES} + q_{(i=1),t}^{Grid} - q_{i,t}^L = \sum_{k\in B_i^D} Q_{ik,t} - \left(Q_{ji,t} - x_{ij}\ell_{ji,t}\right) \quad (10b)$$

$$v_{i,t} = v_{k,t} + 2\left(r_{ik}P_{ik,t} + x_{ik}Q_{ik,t}\right) - \left(r_{ik}^2 + x_{ik}^2\right)\ell_{ik,t} \quad (10c)$$

$$v_{i,t}\ell_{ik,t} = P_{ik,t}^2 + Q_{ik,t}^2 \quad (11)$$

$$-0.6R \leq p_t^{Grid}, q_t^{Grid} \leq R \quad (12)$$

where $i \in B$, $k \in B_i^D$. The subscript ($i = 1$) means the term $p_t^{Grid}$ ($q_t^{Grid}$) only exists in (10a) ((10b)) when $i = 1$.

### C. A Novel Convex Relaxation of DESOS

Before discussing the novel convex relaxation, a technique introduced in [21] is used to eliminate the absolute value sign in (9c). By introducing auxiliary variables $u_{i,t}$ ($i \in B$) which are positive, (9c) can be rewritten as

$$\min \quad f_3 = \sum_t^T \sum_i^N u_{i,t} \quad (13a)$$

$$\text{s.t.} \quad -u_{i,t} \leq v_{i,t} - v_i^{set} \leq u_{i,t}. \quad (13b)$$

In the above DESOS model, the non-convex constraints are (7) and (11). Let $\boldsymbol{x}_{ik,t} = [P_{ik,t}\ Q_{ik,t}\ \ell_{ik,t}\ v_{i,t}]^T$ ($ik \in L$) and $\boldsymbol{y}_{i,t} = [p_{i,t}^{DES}\ q_{i,t}^{DES}\ p_{i,t}^{DESloss}\ v_{i,t}]^T$ ($i \in N_s$), then the convex hulls of (7) and (11) are given as (14) and (15) respectively.

$$\begin{cases}\left\|\boldsymbol{a}^T\boldsymbol{y}_{i,t}\right\|_2 - \boldsymbol{b}^T\boldsymbol{y}_{i,t} \leq 0 \\ \left\|\boldsymbol{a}_i^T\boldsymbol{y}_{i,t}\right\|_2 - \boldsymbol{b}^T\boldsymbol{y}_{i,t} \leq e_i \\ \boldsymbol{c}_i^T\boldsymbol{y}_{i,t} - d_i \leq 0\end{cases} \quad (14)$$

$$\begin{cases}\left\|\boldsymbol{a}^T\boldsymbol{x}_{ik,t}\right\|_2 - \boldsymbol{b}^T\boldsymbol{x}_{ik,t} \leq 0 \\ \boldsymbol{c}_{ik}^T\boldsymbol{x}_{ik,t} - d_{ik} \leq 0\end{cases}, \quad (15)$$

where $\boldsymbol{a}_i = [0\ \sqrt{2r_i^{Batt}}\ 1\ 1]^T$, $\boldsymbol{c}_i = [0\ 0\ \underline{v}_i\bar{v}_i\ r_i^{eq}\left(\bar{S}_i^{DES}\right)^2]^T$, $d_i = r_i^{eq}\left(\bar{S}_i^{DES}\right)^2\left(\bar{v}_i + \underline{v}_i\right)$, and $e_i = r_i^{eq}\left(\bar{S}_i^{DES}\right)^2$. The form of equation (7) is coincidentally similar to (1d). Consequently, following the proof of the theorem in Section II, it is not hard to prove that the set (14) is the convex hull of set (7) within constraints (2b) and (8a). The second convex cut in (14) is required due to the asymmetry of $p_{i,t}^{DES}$ and $q_{i,t}^{DES}$ in (7). As a result, the novel convex relaxation of DESOS is

min    (9a), (9b) or (13)
s.t.     (2), (8), (10), (12), (14) and (15).

### D. Discussion

**The scenario in meshed networks**

Actually, a distribution network is not necessarily radial, especially in the next-generation distribution system [17]. When the considered distribution system is meshed, the BFM given in constraints (1) and (2) are still valid for modeling the AC power flow. However, constraint (1e) which is highly non-convex should be included. Still, convex relaxation (15) is recommended for (1d). For constraint (1e), one can use some existing convex relaxations, like the arctangent envelopes proposed in [5]. However, these convex relaxations may not be necessarily tight due to the high non-convexity of (1e).

**The exactness and sensitive study**

The studies of convex relaxations for the conventional OPF problem is meaningful to many decision-making processes in power systems, like the DESOS problems investigated in this paper, since the AC power flow constraints exist in these problems. However, it has been pointed out in [6] and [7] that, in addition to the tightness, the descent direction of the objective function also determines the exactness of the convex relaxation. Theoretically, the optimal solution of a convex problem with a given (convex) feasible set is determined solely by the objective function. In other words, whether an optimal solution is located in the feasible or infeasible region of the original (non-convex) problem is determined by the chosen objective function. Thus, the descent direction of the objective function should be considered when the exactness of convex relaxations for a given optimization problem is discussed. However, the research efforts in most of the literatures of the convex relaxations for the OPF problem consider only the convex quadratic or linear objective function which represents generation costs. It limits the extendibility of the conclusions obtained in these references to a decision-making process rather than a conventional OPF problem.

In this paper, three objective functions are selected based on the actual operation or planning requirements rather than the requirement of obtaining an exact global optimal solution. The sensitivities of the convex relaxations on the exactness are studied through the cases with the three independent objective functions respectively.

## IV. CASE STUDY

### A. Case Design and Indicators for Performance

In this section, the proposed convex relaxation of DESOS, which is called CH relaxation, is compared with three existing and representative convex relaxations, i.e. the SOCP, basic SDP (BSDP), and enhanced SDP (ESDP) relaxations, through numerical case studies. The SOCP relaxation presented in [4] is used in the comparisons. For the definitions of the BSDP and ESDP relaxations, please refer to reference [7]. Being different from what was defined in [7], the BSDP relaxation considered in this section includes the convex quadratic constraints (2c) and (9a). The ESDP relaxation is constructed by adding the valid linear equalities and the semidefinite inequalities (please refer to (10c) and (12b) in [7] respectively) to the BSDP relaxation.

The optimal objective value (OOV) of a solution is a classical index that has been widely adopted in literature to quantify the tightness of the convex relaxations. Recovering an optimal solution of the convex relaxation to the original problem is another concern. Therefore, in addition to OOV, another indicator is used to capture the feasibility of an obtained solution. For the SOCP and CH relaxations, this indicator is given as

$$e_1^{\max} = \max_{i,ik,t} \left( v_{i,t} \ell_{ik,t} - P_{ik,t}^2 - Q_{ik,t}^2 \right)$$
$$e_2^{\max} = \max_{i,t} \left( p_{i,t}^{DESloss} v_{i,t} - r_i^{eq} \left( p_{i,t}^{DES} \right)^2 - r_i^{Cvt} \left( q_{i,t}^{DES} \right)^2 \right).$$

For the SDP-based relaxations, it is

$$e^{\max} = \max_{i,j,t} \left( X_{ij,t} - x_{i,t} x_{j,t} \right).$$

The above indicators provide measurements of the maximum errors (ME) between the RHS and LHS of the quadratic equality equations with respect to the obtained optimal solutions. If the MEs of a quadratic equation at the optimal solutions are very small values, the relaxation is said to be exact for this equation. In this paper, a ME which is smaller than 0.001 p.u. is considered small enough to claim exactness. The CPU times for solving the convex relaxations are also compared.

### B. Test-bed systems

The convex relaxations of the DESOS algorithm are tested on the IEEE 13, 37, 123-bus feeders [22] assuming that there is high penetration of PV resources and a 9-bus real-world feeder in Arizona [6] respectively. For the three-phase topologies of the four test systems, please refer to Fig. 10 in [7]. The capacities of both PV and DES units for all the test systems are listed in Table III. The problems are all solved by the solver MOSEK (version 7.1.0.53) [23] through the MATLAB toolbox YALMIP [24]. A computer with a 64-bit Intel i5-3230M dual core CPU at 2.60 GHz and 4 GB of RAM was used to run the test cases.

TABLE III
PV SYSTEM AND DES LOCATION AND CAPACITY

| Test system | PV location (bus #) and capacity | penetration |
|---|---|---|
| 9-bus | 4 (0.85 MW), 9 (0.65 MW) | 30.5% |
| 13-bus | 633 (0.5 MW), 680 (0.2 MW), 684 (0.5 MW) | 36.7% |
| 37-bus | 703 (0.3 MW), 706 (0.3 MW), 708 (0.3 MW), 711 (0.3 MW) | 48.8% |
| 123-bus | 8, 15, 25, 44, 54, 67, 81, 89, 105, 110 (all PV systems have the same size of 0.2 MW) | 57.3% |

| Test system | DES location (bus #) and capacity* | |
|---|---|---|
| 9-bus | 4 (1MVA, 2 MWh), 9 (0.75 MVA, 1.5 MWh) | |
| 13-bus | 684 (0.75 MVA, 2.4 MWh), 692 (0.95 MVA, 3.2 MWh) | |
| 37-bus | 720 (0.36 MVA, 1.2 MWh), 730 (0.36 MVA, 1.2 MWh), 737 (0.36 MVA, 1.2 MWh) | |
| 123-bus | 13 (0.23 MVA, 1 MWh), 23 (0.23 MVA, 1 MWh), 76 (0.23 MVA, 1 MWh), 108 (0.23 MVA, 1 MWh) | |

*The capacities of a DES unit include the MVA/kVA capacity of its converter and the energy capacity of the battery.

24-hour data is used in the real-world feeder cases. Since multi-period demand and PV generation profiles are available only for the real-world feeder, the DES dynamic constraint (8b) is omitted and only a snapshot power flow is considered in each of the IEEE feeder cases. When objective function 1 is chosen, an actual 24-hour profile of nodal price (NP) (as shown in Fig. 4) from the website of the ISO New England is used. The coefficient $c$ in objective function 1 is set to be -30 \$/MWh for the IEEE cases.

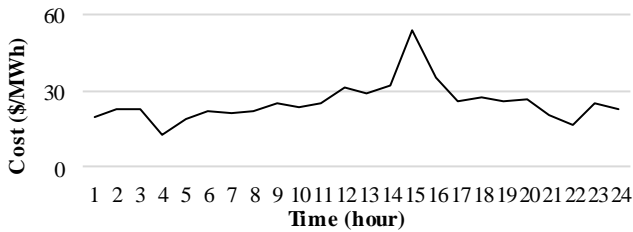

Fig. 4. The 24-hour NP curve for the 9-bus feeder case.

### C. Results and Observations

The results of the cases with objective functions #1-3 are tabulated in Tables IV – VI respectively. The MEs which satisfy the exactness criterion are highlighted in bold. For each case, the exact convex relaxations are highlighted in bold and italic. Based on the numerical results, the following observations are made:

1) The convex relaxations have different performance with different objective functions. When objective function #2 is chosen, the SOCP relaxation is preferred since it is exact for all the cases studied in this paper. However, this does not necessarily mean the SOCP relaxation is exact for all the DESOS cases where objective function 2 is adopted.

2) The ESDP relaxation is tighter than the BSDP relaxation since it is constructed by imposing some more convex constraints on the BSDP relaxation. Similarly, the CH relaxation is tighter than the SOCP relaxation. The SDP-based relaxations neither dominate nor are dominated by the SOCP relaxation when BFM is adopted to describe the AC power flows for DESOS problems. For all the cases studied in this paper, the CH relaxation provides tighter solutions than the ESDP relaxation does.

3) The exact globally optimal solutions for more cases can be achieved by adopting the CH relaxation. In other words, the CH relaxation increases the probability of obtaining the exact globally optimal solutions for actual power systems. This property makes the convex optimization technologies more practical in power systems.

4) The computational burden of CH relaxation is similar to that of the SOCP relaxation which is much smaller than the SDP-based relaxations.

## V. CONCLUSIONS AND FUTURE WORK

The lower bounds on active and reactive power generation [25] and [26] as well as the constraints on voltage magnitudes [27] significantly affect the exactness of convex relaxations for an OPF problem. A case given in [7] shows that the linear state of charging/discharging constraint may also result in inexactness. The descent direction of a chosen objective function also impacts the exactness of the convex relaxation [6] and [7]. As a result, it can be concluded that a large number of factors may determine the exactness of a convex relaxation. It also implies that a convex relaxation which is exact for OPF may not be exact for some other optimization problems that consider the AC power flows as constraints. Studying tighter convex relaxations of the AC power flows is a promising way to mitigate this problem since the AC power flows exist in many optimal decision-making processes in power system.

TABLE IV
RESULTS OF CASE STUDY WITH OBJECTIVE FUNCTION #1

| System | Convex relaxation | OOV ($) | ME#1 | ME#2 | CPU time (s) |
|---|---|---|---|---|---|
| 9-bus real-world feeder | SOCP | 1188.2 | 0.0027 | **1.97E-07** | 0.46 |
| | BSDP | 1170.0 | 8.4345 | | 3.78 |
| | ESDP | 1194.4 | 6.6004 | | 5.77 |
| | *CH* | 1225.3 | **8.10E-04** | **2.10E-07** | 0.50 |
| IEEE 13-bus feeder | SOCP | -89.997 | 8.6353 | 0.5353 | 0.48 |
| | BSDP | -90.000 | 13.9317 | | 3.23 |
| | ESDP | -64.068 | 10.0041 | | 5.08 |
| | *CH* | -34.692 | **4.20E-04** | **6.60E-05** | 0.52 |
| IEEE 37-bus feeder | SOCP | -106.67 | 9.4422 | 0.2815 | 0.49 |
| | BSDP | -106.70 | 40.1427 | | 7.61 |
| | ESDP | -100.50 | 7.7433 | | 17.41 |
| | *CH* | -75.66 | 3.8931 | **7.40E-07** | 0.52 |
| IEEE 123-bus feeder | SOCP | -92.688 | 9.7179 | 0.1344 | 0.56 |
| | BSDP | -95.904 | 27.5926 | | 127.53 |
| | ESDP | -85.806 | 8.4054 | | 203.46 |
| | *CH* | -75.972 | 6.4093 | **1.40E-07** | 0.55 |

TABLE V
RESULTS OF CASE STUDY WITH OBJECTIVE FUNCTION #2

| System | Convex relaxation | OOV (p.u.) | ME#1 | ME#2 | CPU time (s) |
|---|---|---|---|---|---|
| 9-bus real-world feeder | *SOCP* | 0.3514 | **8.00E-04** | **1.10E-07** | 0.58 |
| | BSDP | 0.0000 | 6.5529 | | 3.25 |
| | ESDP | 0.0000 | 3.9984 | | 5.08 |
| | CH | 0.3514 | **8.00E-04** | **1.00E-07** | 0.55 |
| IEEE 13-bus feeder | *SOCP* | 0.0977 | **5.60E-04** | **2.30E-07** | 0.50 |
| | BSDP | 0.0000 | 8.4821 | | 2.30 |
| | ESDP | 0.0000 | 3.5527 | | 4.59 |
| | CH | 0.0995 | **4.20E-04** | **3.90E-07** | 0.56 |
| IEEE 37-bus feeder | *SOCP* | 0.0058 | **5.80E-05** | **2.10E-07** | 0.53 |
| | BSDP | 0.0000 | 52.0979 | | 7.19 |
| | ESDP | 0.0000 | 19.324 | | 12.92 |
| | CH | 0.0058 | **4.90E-05** | **1.90E-07** | 0.48 |
| IEEE 123-bus feeder | *SOCP* | 0.0190 | **2.24E-04** | **1.30E-07** | 0.41 |
| | BSDP | 0.0000 | 40.5173 | | 122.85 |
| | ESDP | 0.0000 | 11.7312 | | 189.07 |
| | CH | 0.0190 | **2.22E-04** | **1.20E-07** | 0.45 |

TABLE VI
RESULTS OF CASE STUDY WITH OBJECTIVE FUNCTION #3

| System | Convex relaxation | OOV (p.u.) | ME#1 | ME#2 | CPU time (s) |
|---|---|---|---|---|---|
| 9-bus real-world feeder | SOCP | 4.1190 | 9.2825 | 0.1701 | 0.52 |
| | BSDP | 4.1293 | 90.3270 | | 4.09 |
| | ESDP | 5.4114 | 76.6040 | | 8.38 |
| | *CH* | 6.8478 | 5.0053 | **1.60E-06** | 0.56 |
| IEEE 13-bus feeder | SOCP | 0.2411 | 8.9359 | 0.5387 | 0.42 |
| | BSDP | 0.2007 | 7.5327 | | 3.19 |
| | ESDP | 0.3413 | 4.1100 | | 5.81 |
| | *CH* | 0.4431 | 1.0888 | **6.30E-06** | 0.39 |
| IEEE 37-bus feeder | SOCP | 0.7056 | 9.2563 | 0.2842 | 0.58 |
| | BSDP | 0.7096 | 25.8243 | | 8.54 |
| | ESDP | 0.7645 | 9.9924 | | 19.84 |
| | *CH* | 0.9777 | 6.0450 | **3.60E-07** | 0.52 |
| IEEE 123-bus feeder | SOCP | 1.4053 | 9.7172 | 0.1330 | 0.44 |
| | BSDP | 1.3978 | 28.7781 | | 135.67 |
| | ESDP | 1.4326 | 9.3852 | | 252.25 |
| | *CH* | 1.4875 | 6.0546 | **1.80E-07** | 0.46 |

The CH relaxation is based on the convex hull of the quadratic non-convex equations in the *DistFlow* model. It can provide a very tight relaxation of AC power flows in radial networks. The case study shows that the CH relaxation works effectively on the DESOS problems considering different objective functions. The CH relaxation provides exact globally optimal solutions for more DESOS cases, which makes the convex optimization more practical in power systems. As a primary future work, the convex hull of the whole *DistFlow* model rather than just the quadratic equation will be studied.

## VI. APPENDIX

Each equality in (1d) is a quadratic polynomial of four variables. It is not possible to visualize its feasible set $\Omega_0$ in the $(P_{ik}, Q_{ik}, \ell_{ik}, v_i)$-space which is a 4-dimensional space. To provide a geometric understanding of the convex hull of $\Omega_0$. Its feasible set is projected to all the 3-dimensional sub-spaces. Note, for example, the projection of $\Omega_0$ onto the $(P_{ik}, Q_{ik}, \ell_{ik})$-space is denoted as $\Omega_{PQ\ell}$.

**Projection in the $(P_{ik}, Q_{ik}, \ell_{ik})$-space**

Considering $v_i$ as a parameter whose value changes from $\underline{v}_i$ to $\overline{v}_i$, equality (1d) can be rewritten as

$$\ell_{ik} = g_1(P_{ik}, Q_{ik}) = (P_{ik}^2 + Q_{ik}^2)/v_i, \quad (16a)$$

which is subjected to (2). The feasible set of (15a), $\Omega_{PQ\ell}$, is sketched in Fig. 5.

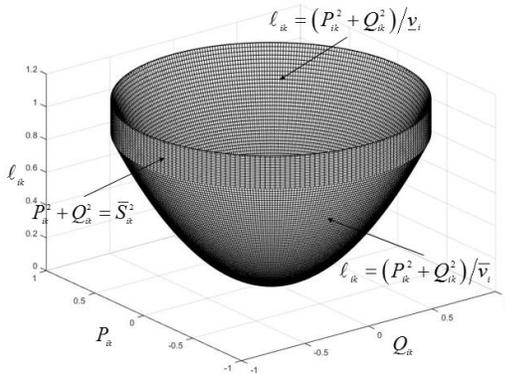

Fig. 5 $\Omega_{PQ\ell}$, the projection of $\Omega_0$ onto the $(P_{ik}, Q_{ik}, \ell_{ik})$-space.

By observation, $CH(\Omega_{PQ\ell})$ can be formulated as

$$P_{ik}^2 + Q_{ik}^2 \leq \overline{v}_i \ell_{ik} \text{ and } \ell_{ik} \leq \overline{\ell}_{ik}. \quad (16b)$$

**Projection in the $(P_{ik}, Q_{ik}, v_i)$-space**

Considering $\ell_{ik}$ as a parameter whose value changes from 0 to $\overline{\ell}_{ik}$, equality (1d) can be rewritten as

$$v_i = g_2(P_{ik}, Q_{ik}) = (P_{ik}^2 + Q_{ik}^2)/\ell_{ik}, \quad (17a)$$

which is subjected to (2). The feasible set of (16a), $\Omega_{PQv}$, is sketched in Fig. 6.

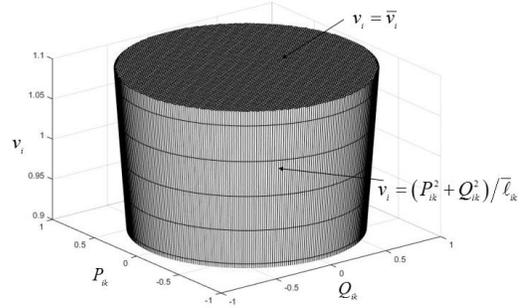

Fig. 6 $\Omega_{PQv}$, the projection of $\Omega_0$ onto the $(P_{ik}, Q_{ik}, v_{ik})$-space.

By observation, $CH(\Omega_{PQv})$ can be formulated as

$$P_{ik}^2 + Q_{ik}^2 \leq \overline{\ell}_{ik} v_i \text{ and } \underline{v}_i \leq v_i \leq \overline{v}_i. \quad (17b)$$

**Projection in the $(P_{ik}, \ell_{ik}, v_i)/(Q_{ik}, \ell_{ik}, v_i)$-space**

In the variable spaces, the positions of $P_{ik}$ and $Q_{ik}$ are completely symmetrical. As a result, the formulation of $CH(\Omega_{P\ell v})$ can be directly applied to obtain $CH(\Omega_{Q\ell v})$ by replacing $P_{ik}$ with $Q_{ik}$. Considering $Q_{ik}$ as a parameter, equality (1d) can be rewritten as

$$P_{ik} = g_3(\ell_{ik}, v_i) = \pm\sqrt{v_i \ell_{ik} - Q_{ik}^2}, \quad (18a)$$

which is subjected to (2). The feasible set of (17a), $\Omega_{P\ell v}$, is sketched in Fig. 7.

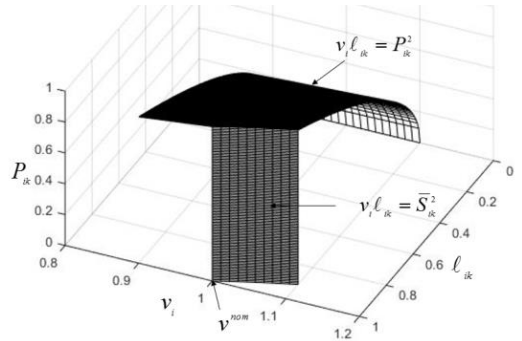

Fig. 7 $\Omega_{P\ell v}$. $\Omega_{P\ell v}$ is the space that is enclosed by the bounding planes described by (2), and the two surfaces shown in the figure.

By observation, inside space (2), $CH(\Omega_{P\ell v})$ can be formulated as

$$\begin{cases} \left\| \begin{array}{c} \sqrt{2} P_{ik} \\ v_i - \ell_{ik} \end{array} \right\|_2 \leq v_i + \ell_{ik} \\ \overline{v}_i \ell_{ik} + \overline{\ell}_{ik} v_i \leq \overline{\ell}_{ik} (\overline{v}_i + v^{nom}) \end{cases} \quad (18b)$$

Consequently, inside space (2), $CH(\Omega_{Q\ell v})$ can be expressed as

$$\begin{cases} \left\| \begin{array}{c} \sqrt{2} Q_{ik} \\ v_i - \ell_{ik} \end{array} \right\|_2 \leq v_i + \ell_{ik} \\ \overline{v}_i \ell_{ik} + \overline{\ell}_{ik} v_i \leq \overline{\ell}_{ik} (\overline{v}_i + v^{nom}) \end{cases}. \quad (18c)$$

It is easy to observe that

$$\begin{cases} (\Omega_2)_{PQ\ell} = CONV(\Omega_{PQ\ell}) \\ (\Omega_2)_{PQv} = CONV(\Omega_{PQv}) \\ (\Omega_2)_{\ell vP} = CONV(\Omega_{\ell vP}) \\ (\Omega_2)_{\ell vQ} = CONV(\Omega_{\ell vQ}) \end{cases},$$

where $(\Omega_2)_{PQ\ell}$ means the projection of $\Omega_2$ onto the ($P_{ik}$ $Q_{ik}$ $v_i$)-space.

The projection of a convex set onto a subspace should also be convex. However, convex projections on subspace do not sufficiently mean that the original set is convex. Thus, the above geometric study provides a necessary but not sufficient validation of the Theorem in Section II. Nevertheless, by providing a visual interpretation, it may be of help to readers in understanding the Theorem.